\documentclass{amsart}

\usepackage{etex}
\usepackage{amsmath, amssymb}
\usepackage{array}
\usepackage[frame,cmtip,arrow,matrix,line,graph,curve]{xy}
\usepackage{graphpap, color, paralist, pstricks}
\usepackage[mathscr]{eucal}
\usepackage[pdftex]{graphicx}
\usepackage[pdftex,colorlinks,backref=page,citecolor=blue]{hyperref}
\usepackage{pifont}
\usepackage{cancel}
\usepackage{mathtools}
\usepackage{tikz}

\usepackage{bbm}
\usepackage{dsfont}

\newcommand{\pr}{\mathbb{P}}

\newcommand{\beq}[1]{\begin{equation}\label{#1}}
\newcommand{\enq}[0]{\end{equation}}

\newcommand{\mn}[0]{\medskip\noindent}
\newcommand{\nin}[0]{\noindent}

\newcommand{\sub}[0]{\subseteq}
\newcommand{\sm}[0]{\setminus}
\renewcommand{\dots}[0]{,\ldots,}

\newcommand{\ov}[0]{\overline}

\newcommand{\A}[0]{{\mathcal A}}
\newcommand{\B}[0]{{\mathcal B}}
\newcommand{\cee}[0]{{\mathcal C}}
\newcommand{\D}[0]{{\mathcal D}}
\newcommand{\eee}[0]{{\mathcal E}}

\newcommand{\Lra}[0]{\Leftrightarrow}

\newcommand{\0}[0]{\emptyset}

\newcommand{\ga}[0]{\alpha }

\newcommand{\go}[0]{\omega}
\newcommand{\gO}[0]{\Omega}

\newcommand{\gs}[0]{\sigma}

\newcommand{\prh}[1][]{\pr_h}
\newcommand{\oH}[1][]{\ov{H}}
\newcommand{\oB}[1][]{\ov{B}}

\newtheorem{thm}{Theorem}

\newtheorem{cor}[thm]{Corollary}

\newtheorem{conj}[thm]{Conjecture}
\theoremstyle{definition}

\newtheorem{question}[thm]{Question}

\newcommand{\mbone}[0]{\mathbbm{1}}

\begin{document}

\title{A note on positive association}

\author{Jeff Kahn}
\thanks{Department of Mathematics, Rutgers University}
\thanks{Supported by NSF Grant DMS1954035 and a Simons Fellowship}
\email{jkahn@math.rutgers.edu}
\address{Department of Mathematics, Rutgers University \\
Hill Center for the Mathematical Sciences \\
110 Frelinghuysen Rd.\\
Piscataway, NJ 08854-8019, USA}

\begin{abstract}
We show 
that if $\A,\B,\cee$ are increasing subsets of $\gO:=\{0,1\}^n$
with $\A\neq\0$, then with respect to any product probability measure on $\gO$,
\[
\mbox{if each of the pairs $\{\A\cap\B,\cee\}$, $\{\A\cap \cee,\B\}$ is independent,
then $ \B$ and $\cee$ are independent.}
\]
This implies an answer to a motivating question of J.\ Steif,
and is related to a basic, still open variant of that question, and to a well-known 
conjecture of S.\ Sahi.
\end{abstract}

\maketitle

\section{Introduction}\label{Intro}
This began with a question of Jeff Steif that I first heard in
conversation with him and Rob van den Berg many years ago.  
(The participants have been unable to agree as to the decade in which this
conversation took place, or the continent that hosted it, but it was 
no later than 2002.)
The aims of this note are to:

\mn
(i) prove a small result (Theorem~\ref{TABC}) that implies an answer to 
Steif's question;

\mn
(ii) point out that a variant (Question~\ref{steifQ}), which seems as basic a question
as one could ask about positive association, remains open; and

\mn
(iii) observe a curious connection to a well-known conjecture of Siddhartha Sahi.

\mn
\emph{Background.}
Throughout this discussion,
$X_1\dots X_n$ are Bernoulli random variables and $\mu$ is the 
law of $(X_1\dots X_n)$;
so $\mu$ is a probability measure on 
$\gO:=\{0,1\}^n\equiv 2^{[n]}$ (with the natural identification of sets with their indicators).
Recall that events $\A,\B$ (in \emph{any} probability space)
are \emph{positively
correlated}
if $\pr(\A\B)\geq\pr(\A)\pr(\B)$.
The law of 
$X_{1},...,X_{n}$
is \emph{positively associated} (PA),
or has {\em positive association}, 
if 
\[
\mbox{\emph{any two events both increasing in the $X_{i}$'s are positively correlated.}}
\]

The seminal result here is {\em Harris' Inequality} \cite{Harris}, which says
\emph{product measures are PA.}
(To be precise, this is what's given by Harris' argument;
the statement in \cite{Harris} is less general.)
Harris' Inequality for uniform measure 
was rediscovered in \cite{Kleitman}, and in combinatorial circles is sometimes
called {\em Kleitman's Lemma}.
The most useful extension of Harris (discovered still later but still independently) is the
\emph{FKG Inequality}
of Fortuin, Kasteleyn, and Ginibre
\cite{FKG},
which says $\mu$ is PA whenever
\beq{PLC}
\mu(A)\mu(B)\leq\mu(A\cap B)\mu(A\cup B)
\,\,\forall A,B\in\Omega.
\enq
A $\mu$ satisfying (\ref{PLC}) 
(the ``positive lattice condition")
is an {\em FKG measure}.
See e.g.\ 
\cite{GrimmettRC}, \cite{Liggett} for some indication of the role of positive association
in probability, and \cite[Ch.\ 6]{Alon-Spencer} for a quick hint on the combinatorial side.

\emph{For the rest of the paper, script capitals ($\A,\B,\ldots$) are 
nonempty increasing events in $\gO$.}
We use $\A\B$ for $\A\cap \B$, and denote independence of $\A$ and $\B$ 
by $\A|\B$ and dependence by $\A\sim \B$.

\mn
\emph{Underlying independents.}
One way to prove PA for $\mu$ 
is to realize the $X_i$'s as increasing
functions of
{\em independent} Bernoullis $Y_1\dots Y_m $ and invoke Harris; 
more generally, $\mu$ is PA if it is a limit of measures obtained in
this way.
Say $\mu$ is {\em FUI}
(for {\em finitely many underlying independents)}
in the first case, and {\em UI} in the second.
This is not as restrictive as it sounds, since
\emph{FKG measures are FUI},\footnote{Since it seems hard to find a reference for this,
I include here a sketch of a proof that was shown to me by Rob van den Berg,
which he believes is (implicitly) well known in the probability community:
(a)  Assuming the law, $\mu$, of $(X_1\dots X_n)$ is FKG, let $Z_1\dots Z_n$ be independent,
each uniform from $[0,1]$, and for $i=1\dots n$, if $X_j=\go_j$ for $j<i$,
let $X_i=1$ iff $Z_i>1- \mu(X_i=0|X_j=\go_j\forall j<i)$.
This is easily seen to return $\mu$ as the law of $(X_1\dots X_n)$, and it's not hard to
see, using \eqref{PLC}, that the $X_i$'s are nondecreasing in the $Z_j$'s. 
(b) For each $i$, the procedure in (a) depends on a finite number of events 
$A(i,j):=\{Z_i>\ga_{i,j}\}$, and it's easy to realize the indicators $\mbone_{A(i,j)}$
($j\in [m]$, say) as nondecreasing functions of independent (nonidentical) Bernoullis $Y_{i,j}$ 
($j\in [m]$).}
\footnote{For completeness we note that FUI does not imply FKG.
E.g.\ with $Y_1,Y_2,Y_3$ independent Bernoullis,
let $X_1=Y_1Y_2$, $X_2=Y_1Y_3$, and $X_3=Y_2Y_3$, and notice that 
the law of the $X_i$'s assigns probability zero to strings of weight 2,
so trivially violates \eqref{PLC}.}
and in fact Steif's original question was 
\beq{Steif}
\mbox{\emph{are all PA measures FUI?}}
\enq
As we will see shortly, the answer is no; but, remarkably,
we can't (as far as I know) rule out a slightly weaker possibility:
\begin{question}\label{steifQ}
\emph{Are all PA measures UI?}
\end{question}
\nin
Of course one hopes the answer to this very basic question is again no---that is, 
positive association is more than Harris' Inequality---but it seems 
surprisingly hard to say anything about the law of a UI $\mu$ that uses
more than positive association.  In contrast, the following statement,
the technical content of the present note, does 
manage to
distinguish FUI from PA, 
and to imply the promised 
negative answer to \eqref{Steif}.

\begin{thm}\label{TABC}
For any FUI $\mu$ and
(increasing) $\A,\B,\cee$,
\beq{ABC}
\mbox{if $\A\B|\cee$ and $\A\cee |\B$ (and $\mu(\A)\neq 0$),
then  $ \B|\cee$.}
\enq
\end{thm}
\nin

For the connection
to \eqref{Steif}, we recall a 
beautiful result of Doyle, Fishburn and Shepp \cite{DFS}:
\begin{thm}\label{TDFS}
For a uniform permutation $\gs$ of $[n]$, the law, $\mu_n$, of the set of fixed points
of $\gs$ (that is, of $(X_1\dots X_n)$, where $X_i=\mbone_{\{\gs(i)=i\}}$) is PA.
\end{thm}

\nin
\begin{cor}
The answer to \eqref{Steif} is negative.
\end{cor}
\nin
\emph{Proof.}  This follows from Theorem~\ref{TABC} and the observation that
$\mu=\mu_3$ violates \eqref{ABC}:
$\mu$ assigns weight $1/3$ to $\0$
and $1/6$ to each of $\{1\}$, $\{2\}$, $\{3\}$, $\{1,2,3\}$
(and 0 to pairs); so, with $\A_{i,j}= \{\mbox{$\gs$ fixes at least one of $i,j$}\}$,
$\A=\A_{1,2}$, $\B=\A_{1,3}$ and $\cee=\A_{2,3}$, we have 
$\mu(\A)=\mu(\B)=\mu(\cee) =1/2$, $\mu(\A\B) = \cdots = 1/3$, and $\mu(\A\B\cee)=1/6$,
whence the hypotheses of \eqref{ABC} hold but the conclusion does not.\qed

\mn
\emph{Aside}.  As its discoverers emphasize, the argument of \cite{DFS} is a quite
painful case analysis.  Shouldn't there be a nicer, more enlightening
proof of such an elegant result?

\mn
\emph{Sahi's Conjecture.}
This fascinating (infuriating) conjecture \cite{Sahi}
proposes an extension of Harris' Inequality to $k>2$ events;
we state just the case $k=3$, which has to date proved thoroughly intractable
and seems not unlikely to capture the full 
difficulty of the problem.

\begin{conj}\label{richards}
For a product measure $\mu$ and increasing
$\A,\B,\cee\sub \gO$,
\beq{sahi}
2\mu(\A\B\cee)-[\mu(\A\B)\mu(\cee)+\mu(\A\cee)\mu(\B)+\mu(\B\cee)\mu(\A)]
+\mu(\A)\mu(\B)\mu(\cee)\geq 0.
\enq
\end{conj}
\nin
Note $\mu(\A)=1$ recovers Harris.
The conjecture is stated in \cite{Sahi} for FKG measures, but this is no more general
since FKG measures are FUI.
For $k\in\{3,4,5\}$, Sahi's Conjecture was originally stated---as a theorem, but with an incorrect 
proof---by Richards \cite{Richards}; he also suggested the possibility of similar 
inequalities for larger $k$, but without proposed candidates for the coefficients.
Progress on the conjecture has been limited (see \cite{Lieb-Sahi} for the state of the art
and \cite{Sahi2,Sahi3} for related results), surely a poor reflection of the
effort expended on it.

For present purposes the point of all this is that \eqref{sahi}, with Harris,
implies \eqref{ABC} (since under the hypotheses of \eqref{ABC}, \eqref{sahi}
becomes $\mu(\B\cee)\leq \mu(\B)\mu(\cee)$); so a positive answer to 
Question~\ref{steifQ}, even just for $\mu_3$, would say Sahi's Conjecture---which
of course also implies \eqref{sahi} when $\mu$ is UI---is false.
Conversely, Theorem~\ref{TABC} may be considered a tiny step \emph{toward} 
Sahi's Conjecture.

\mn

In Section~\ref{Proofs} we will give two proofs of Theorem~\ref{TABC}.
The first of these is very short and rather \emph{ad hoc}.
The second is longer (not long) but feels more systematic, and is 
included here in the hope that it might be more susceptible to improvement.

\section{Proofs}\label{Proofs}

As usual, $\min(\A)$ is the set of minimal elements of $\A$. 
One says $i\in [n]$ \emph{affects} $\A$ if 
$A\cup \{i\} \in \A$ for some $A\not\in \A$
(and $I$ \emph{affects} $\A$ if some $i\in I$ does).  We use $Z(\A)$ for the set of $i\in [n]$ that affect $\A$, 
noting that
\[
Z(\A) =\cup \{A:A\in \min(\A)\}.
\]

The basis for both our proofs of Theorem~\ref{TABC}, an
immediate consequence of Harris' argument, is 
\beq{basic}
\A|\B ~\Lra ~Z(\A)\cap Z(\B) =\0.
\enq
(Equivalently, $\A$ and $\B$ are independent iff $A\cap B=\0$ whenever 
$A\in \min(\A)$ and $B\in \min(\B)$.)

\mn

For each of the following arguments 
we assume $\B \sim\cee$ (and $\mu(\A)\neq 0$) and want to show
\beq{at least one}
\mbox{at least one of $\A\B\sim\cee$,
$\A\cee\sim\B$ holds;}
\enq
so we assume \eqref{at least one} fails and aim for a contradiction.

\mn
\emph{\textbf{First proof.}}
(A Venn diagram may be helpful here.)
Notice to begin that, for any increasing $\D$ and $\eee$,
\[
\mbox{$\min(\D\eee) $ is the set of minimal elements of 
$ \{D\cup E:D\in \min(\D),E\in \min(\eee)\}$.}
\]
We now use $A$ (possibly subscripted) for members of $\min(\A)$ and so on.

Recalling that we assume $\B\sim \cee$, 
choose $B$ and $C$ with $B\cap C\neq \0$ and $B\cup C$ minimal subject to this,
and observe (with justification below) that 
\[
B\cup C\in \min(\B\cee).
\]
\emph{Proof.}
Suppose instead that $B\cup C\supsetneq B_0\cup C_0$, 
say with $B_0\neq B$.
Then $B_0\cap (C\sm B)\neq\0$ (since otherwise $B_0\subsetneq B$).
On the other hand, $C_0\cap B\neq \0$ (else $C\neq C_0\sub C$) 
implies $C_0\supseteq C\sm B$ (else $B\cup C_0\subsetneq B\cup C$ contradicts
our choice of $(B,C)$).
But then $B_0\cap C_0\neq \0$, which is again a contradiction.\qed

\mn

Choose $A$ with $A\sm (B\cup C)$ minimal.
Since (we assume)
$\A\B |\cee$, there must be some $A_1\cup B_1\in \min(\A\B)$ 
with $A_1\cup B_1\sub (A\cup B)\sm C$;
in particular $B_1\sub (A\cup B)\sm C$, implying $\0\neq B_1\sm B\sub A\sm (B\cup C)$.
But then, since $\A\cee |\B$, there is $A_2\cup C_2\sub (A\cup C)\sm B_1$; which
contradicts our choice of $A$ since
$A_2\sm (B\cup C)\sub A\sm (B\cup C\cup  B_1) \subsetneq A\sm (B\cup C)$.

~
\hfill $\blacksquare$

\mn
\emph{\textbf{Second proof.}}
Let $I=Z(\B)\cap Z(\cee)$ ($\neq \0)$, 
$J = Z(\B)\sm Z(\cee)$, $K =Z(\cee)\sm Z(\B)$, and $L=[n]\sm (I\cup J\cup K)$.

We are (again) assuming \eqref{at least one} fails, so in particular,
\beq{Idoesnt}
\mbox{$I$ doesn't affect either of $\A\B$, $\A\cee$.}
\enq
Of course we may also assume 
\beq{AsubBcup}
\A\sub \B\cup\cee,
\enq 
since replacing $\A$ by $\A\cap (\B\cup \cee)$ has no effect on $\A\B$, $\A\cee$.

\mn
\emph{Observation 1.}
$I\cap Z(\A)=\0$.

\mn
\emph{Proof.}
Suppose $i\in I$, $A\not\in \A$, $A':=A\cup \{i\}\in \A$, and (w.l.o.g.; see \eqref{AsubBcup})
$A'\in \B$.  Then $\A\B\cap \{A,A'\}=\{A'\}$, implying $i\in Z(\A\B)$ and
contradicting \eqref{Idoesnt}.
\qed

\mn
\emph{Observation 2.}
If $X\sub J$, $X\not\in\B$ and $X\cup I\in \B$, then $X\cup K\cup L\not\in \A$
(so $X\cup Y\not\in \A$ $~\forall Y\sub K\cup L$).

\mn
(And similarly with $\B$ replaced by $\cee$ and the roles of $J$ and $K$ interchanged.)

\mn
\emph{Proof.}
Otherwise $X\cup K\cup L\cup I\in \A\B$ and $X\cup K\cup L\not\in \A\B$ 
(since $X\not\in \B$ and $(K\cup L)\cap Z(\B)=\0$), contradicting \eqref{Idoesnt}.
\qed

\mn

Notice that 
$I\sub Z(\B)= I\cup J$
implies that there is some $X\sub J$ with
$\B\cap \{X,X\cup I\}=\{X\cup I\}$.
Let $X$ be 
of this type and, similarly,
let $Y\sub K$ 
satisfy $Y\not\in \cee$ and $Y\cup I\in \cee$.  

Let $A$ be minimal in $\A$ with $A\supseteq X\cup Y$.  Then
$A\cap I=\0$ (by Observation 1),
so by Observation 2,
\beq{AJX}
\mbox{$A\cap J\supsetneq X~$ and $~A\cap K\supsetneq Y$.}
\enq
(We just need one of these.) 
Now
\[
A\cup I\in \A\B\cee
\]
(since $X\cup I\in \B$ and $Y\cup I\in \cee$; we just need $A\cup I\in \A\cee$), 
while minimality of $A$ and Observation 1 give 
\[
(A\cup I)\sm \{j \}\not\in \A
\]
for any $j\in (A\cap J)\sm X$ (and \eqref{AJX} says there \emph{is} such a $j$).
Thus any such $j$ is in $Z(\A\cee)\cap Z(\B)$, so that, contrary to assumption, 
\eqref{at least one} \emph{does} hold.

~
\hfill $\blacksquare$

\mn
\textbf{Acknowledgment.}
This paper owes its existence to the interest of Rob van den Berg and Jeff Steif.\linebreak
Theorem~\ref{TABC} dates to the time of the meeting between the three of us mentioned 
earlier (the second proof in Section~\ref{Proofs} is newer); 
but it was only recently that they convinced me that the result---and the discussion in 
Section~\ref{Intro}---ought to be published, and I 
am greatly indebted to
them for their encouragement, 
and their many helpful comments on the manuscript.

\end{document}